\documentclass[12pt]{article}

\usepackage{amsmath,amssymb,enumerate,bbm}
\newcommand{\upl}{+}
\newenvironment{enumeratealpha}{\begin{enumerate}[a{\textup{)}}] }{\end{enumerate}}
\newenvironment{enumeratenumeric}{\begin{enumerate}[1.] }{\end{enumerate}}
\newenvironment{proof}{\noindent\textbf{Proof\ }}{\hspace*{\fill}$\Box$\medskip}
\newtheorem{theorem}{Theorem}

\newtheorem{lemma}{Lemma}

\begin{document}

\title{Explicit Ramsey graphs and Erd\"os distance problem over finite Euclidean and non-Euclidean spaces}\author{Le Anh Vinh\\
Mathematics Department\\
Harvard University\\
Cambridge, MA 02138, US\\
vinh@math.harvard.edu}\maketitle

\begin{abstract}
 We study the Erd\"os distance problem over finite Euclidean and non-Euclidean spaces. Our main tools are graphs associated to finite Euclidean and non-Euclidean spaces that are considered in Bannai-Shimabukuro-Tanaka (2004, 2007). These graphs are shown to be asymptotically Ramanujan graphs. The advantage of using these graphs is twofold. First, we can derive new lower bounds on the Erd\"os distance problems with explicit constants. Second, we can construct many explicit tough Ramsey graphs $R(3,k)$. 
\end{abstract}

\section{Introduction}

Let $\mathbbm{\mathbbm{F}}_q$ denote the finite field with $q$ elements where
$q \gg 1$ is an odd prime power. Let $E \subset \mathbbm{F}_q^d$, $d \geqslant
2$. Then the analog of the classical Erd\"os distance problem is to determine
the smallest possible cardinality of the set
\[ \Delta (E) =\{|x - y|^2 = (x_1 - y_1)^2 + \ldots + (x_d - y_d)^2 : x, y \in
   E\}, \]
viewed as a subset of $\mathbbm{F}_q$. Suppose that $- 1$ is a square in
$\mathbbm{F}_q$, then using spheres of radius $0$, there exists a set of
cardinality precisely $q^{d / 2}$ such that $\Delta (E) =\{0\}$. Thus, we only
consider the set $E \subset \mathbbm{F}_q^d$ of cardinality $C q^{\frac{q}{2}
+ \varepsilon}$ for some constant $C$. Bourgain, Katz and Tao (\cite{bourgain-katz-tao}) showed, using
intricate incidence geometry, that for every $\varepsilon > 0$, there exists
$\delta > 0$, such that if $E \in \mathbbm{F}_q^2$ and $|E| \leqslant C_{\epsilon} q^{2 -
\epsilon}$, then $ |\Delta (E)| \geqslant C_{\delta} q^{\frac{1}{2} +
\delta}$ for some constants $C_{\epsilon}, C_{\delta}$. The relationship
between $\varepsilon$ and $\delta$ in their argument is difficult to
determine. Going up to higher dimension using arguments of Bourgain, Katz and
Tao is quite subtle. Iosevich and Rudnev (\cite{iosevich-rudnev}) establish the following results
using Fourier analytic methods.

\begin{theorem}\label{ir1} (\cite{iosevich-rudnev})
  Let $E \subset \mathbbm{F}_q^d$ such that $|E| \gtrsim C q^{d / 2}$ for $C$
  sufficient large. Then
  \begin{equation}
    |\Delta (E)| \gtrsim \min \left\{ q, \frac{|E|}{q^{\frac{d - 1}{2}}}
    \right\} .
  \end{equation}
\end{theorem}

By modifying the proof of Theorem \ref{ir1} slightly, they obtain the following stronger
conclusion.

\begin{theorem} (\cite{iosevich-rudnev})
  Let $E \subset \mathbbm{F}_q^d$ such that $|E| \geqslant C q^{\frac{d +
  1}{2}}$ for sufficient large constant $C$. Then for every $t \in
  \mathbbm{F}_q$ there exist $x, y \in E$ such that $|x - y|^2 = t$. In other
  words, $|\Delta (E)| = q$.
\end{theorem}

It is, however, more natural to define the analogues of Euclidean
graphs for each non-degenerate quadratic from on $V= \mathbbm{F}_q^d$, $d \geqslant
2$. Let $Q$ be a non-degenerate quadratic form on $V$. For any $E \subset V$, we define the distance set of $E$ with respect to $Q$:
\[ \Delta_Q (E) =\{Q(x - y) : x, y \in E\}, \]
viewed as a subset of $\mathbbm{F}_q$. Our first result is the following.

\begin{theorem}\label{v1} Let $Q$ be a non-degenerate quadratic from on $\mathbbm{F}_q^d$, $d \geqslant
2$.   Let $E \subset \mathbbm{F}_q^d$ such that $|E| \geqslant 3 q^{\frac{d}{2} +
  \varepsilon}$ for some $\varepsilon > 0$, then
  \begin{equation}
    |\Delta_Q (E)| \geqslant \min \{ \frac{|E|}{3 q^{(d - 1) / 2}}, q\}
  \end{equation}
  for $q \gg 1$.
\end{theorem}

An interesting question is to study the analogous of the Erd\"os distance problem in non-Euclidean spaces. In order to make this paper concise, we will only consider the Erd\"os distance problem in the finite non-Euclidean plane (or so-called the finite upper half plane). In Section 2, we will see how to obtain various finite non-Euclidean spaces from the action of classical Lie groups
on the set of non-isotropic points, lines and hyperplanes. Most of our results in this paper hold in this more general setting. We will address these results in a subsequent paper.

The well-known finite upper half plane constructed in a similar way using an analogue of Poincar\'e's non-Euclidean distance. We follow the construction in \cite{survey}. Let $\mathbbm{F}_q$ be the finite field with $q = p^r$ elements, where $p$ is an odd prime. Suppose $\sigma$ is a generator of the multiplicative group $\mathbbm{F}_q^*$ of nonzero elements in $\mathbbm{F_q}$. The extension $\mathbbm{F}_q \cong \mathbbm{F}_q(\sigma)$ is analogous to $\mathbbm{C}=\mathbbm{R}[i]$. We define the \textit{finite Poincar\'e upper half-plane} as 
\begin{equation} H_q = \{z=x+y\sqrt{\sigma}: x, y \in \mathbbm{F}_q \, \text{and} \, y \neq 0 \}.\end{equation}
Note that ``half-plane'' is something of a misnomer since $y \neq 0$ may not be a good finite analogue of the condition $y > 0$ that defines the usual Poincar\'e upper half-plane in $C$. In fact, $H_q$ is more like a double covering of a finite upper half-plane. We use the familiar notation from complex analysis for $z = x+y \sqrt{\sigma} \in H_q$: $x= Re(z)$, $y= Im(z)$, $\bar{z} = x- y\sqrt{\sigma} = z^q$, $N(z) = $ Norm of $z = z\bar{z} = z^{1+q}$. The \textit{Poincar\'e distance} between $z, w \in H_q$ is
\begin{equation} d(z,w) = \frac{N(z-w)}{Im(z)Im(w)} \in \mathbbm{F}_q.
\end{equation}
This distance is not a metric in the sense of analysis, but it is $GL(2,\mathbbm{F}_q)$-invariant: $d(gz,gw) = d(z,w)$ for all $g \in GL(2,\mathbbm{F}_q)$ and all $z,w \in H_q$. Let $E \subset H_q$. We define the distance set with respect to the Poincar\'e distance:
\[ \Delta_H (E) =\{d(x,y) : x, y \in E\}, \]
viewed as a subset of $\mathbbm{F}_q$. The following result is a non-Euclidean analog of Theorem \ref{v1}.

\begin{theorem}\label{v2} Let $E \subset H_q$ such that $|E| \geqslant 3 q^{\frac{1}{2} +
  \varepsilon}$ for some $\varepsilon > 0$, then
  \begin{equation}
    | \Delta_H (E)| \geqslant \min \{ \frac{|E|}{3 q^{1/2}}, q-1\}
  \end{equation}
  for $q \gg 1$.
\end{theorem}

We also have the Erd\"os problem for two sets. Let $E, F \subset \mathbbm{F}_q^d$, $d \geqslant 2$. Given a non-degenerate quadratic $Q$ form on $\mathbbm{F}_q^d$. We define the set of distances between two sets $E$ and $F$:
\[ \Delta_Q(E,F) = \{Q(x,y) : x \in E, y \in F \}.\] 

We will prove the following analogues of Theorem \ref{v1} for the distance set $\Delta_Q(E,F)$. 

\begin{theorem}\label{v3}
Let $E, F \subset \mathbbm{F}_q^d$ such that $|E||F| \geqslant 9q^{(d-1)+\epsilon}$ for some $\varepsilon > 0$, then
\[ \Delta_Q(E,F) \geqslant \min \{ \frac{\sqrt{|E||F|}}{3q^{(d-1)/2}}, q \}\]
for $q \gg 1$.
\end{theorem}

In finite upper half plane, we define the set of distances between two sets $E, F \subset H_q$:
\[ \Delta_H(E,F) = \{d(x,y) : x \in E, y \in F \},\]
where $d(x,y)$ is the finite Poincar\'e distance between $x$ and $y$. Similarly, we have an analog of Theorem \ref{v2} for the distance set $\Delta_H(E,F)$.

\begin{theorem}\label{v4}
Let $E, F \subset H_q$ such that $|E||F| \geqslant 9q^{1+2\epsilon}$ for some $\varepsilon > 0$, then
\[ \Delta_H(E,F) \geqslant \min \{ \frac{\sqrt{|E||F|}}{3q^{1/2}}, q-1\}\]
for $q \gg 1$.
\end{theorem}

Note that an analog version of Theorem \ref{v3} for $q$ prime was obtained by Iosevich and Koh in \cite{iosevich-koh} (part b of Theorem 3.4) and Shparlinski in \cite{shparlinski}. However, the bound we obtained here is better than ones in these papers. 

The rest of this paper is organized as follows. In Section 2 we construct our main tools to study the Erd\"os problem over finite Euclidean and non-Euclidean spaces, the finite Euclidean and non-Euclidean graphs. Our construction follows one of Bannai, Shimabukuro and Tanaka in \cite{bannai 1, bannai 2}. In Section 3 we establish some useful facts about these finite graphs. One important result is for infinitely many values of $q$, these graphs disprove a conjecture of Chvat\'al and also provide a good lower bound for the Ramsey number $R(3,k)$. We then prove our main results, Theorems \ref{v1}, \ref{v2}, \ref{v3} and \ref{v2}, in Section 4. 

We also call the reader's attention to the fact that the application of the spectral method from graph theory in sum-product estimates and Erd\"os distance problem was independently used by Vu in \cite{van}. 

\section{Finite Euclidean and non-Euclidean Graphs}

In this section, we summarise main results from Bannai-Shimabukuro-Tanaka \cite{bannai 2, bannai 1}. We follow their
constructions of finite Euclidean and non-Euclidean graphs.

Let $Q$ be
a non-degenerate quadratic form on $V$. We define the corresponding bilinear
from on $V$:
\[ \left\langle x, y \right\rangle_Q = Q (x + y) - Q (x) - Q (y) . \]
Let $O (V, Q)$ be the group of all linear transformations on $V$ that fix $Q$
(which is called the orthogonal group associated with the quadratic form $Q$).
The non-degenerate quadratic forms over $\mathbbm{F}_q^n$ are classified as
follows:
\begin{enumeratenumeric}
  \item Suppose that $n = 2 m$. If $q$ odd then there are two inequivalent
  non-degenerate quadratic forms $Q^+_{2 m}$ and $Q^-_{2 m}$:
  \begin{eqnarray*}
    Q^+_{2 m} (x) & = & 2 x_1 x_2 + \ldots + 2 x_{2 m - 1} x_{2 m},\\
    Q^-_{2 m} (x) & = & 2 x_1 x_2 + \ldots + 2 x_{2 m - 3} x_{2 m - 2} + x_{2
    m - 1}^2 - \alpha x_{2 m}^2,
  \end{eqnarray*}
  where $\alpha$ is a non-square element in $\mathbbm{F}_q$. If $q$ even then
  there are also two inequivalent non-degenerate quadratic forms $Q^+$ and
  $Q^-$:
  \begin{eqnarray*}
    Q^{\upl}_{2 m} (x) & = & x_1 x_2 + \ldots + x_{2 m - 1} x_{2 m,}\\
    Q^-_{2 m} (x) & = & x_1 x_2 + \ldots + x_{2 m - 3} x_{2 m - 2} + x_{2 m -
    1}^2 + \beta x_{2 m}^2,
  \end{eqnarray*}
  where $\beta$ is an element in $\mathbbm{F}_q$ such that the polynomial $t^2
  + t + \beta$ is irreducible over $\mathbbm{F}_q$. We write $O^+_{2 m} = O
  (V, Q^+_{2 m})$ and $O^-_{2 m} = O (V, Q^-_{2 m})$.
  
  \item Suppose that $n = 2 m + 1$ is odd. If $q$ is odd, then there are two
  inequivalent non-degenerate quadratic forms $Q_{2 m + 1}$ and $Q'_{2 m +
  1}$:
  \begin{eqnarray*}
    Q_{2 m + 1} (x) & = & 2 x_1 x_2 + \ldots + 2 x_{2 m - 1} x_{2 m} + x_{2 m
    + 1}^2,\\
    Q'_{2 m + 1} (x) & = & 2 x_1 x_2 + \ldots + 2 x_{2 m - 1} x_{2 m} +
    \alpha x_{2 m - 1}^2,
  \end{eqnarray*}
  where $\alpha$ is a non-square element in $\mathbbm{F}_q$. But the groups $O
  (V, Q_{2 m + 1})$ and $O (V, Q'_{2 m + 1})$ are isomorphic. If $q$ is even
  then there exists exactly one inequivalent non-degenerate quadratic form
  $Q_{2 m + 1}$:
  \[ Q_{2 m + 1} (x) = x_1 x_2 + \ldots + x_{2 m - 1} x_{2 m} + x_{2 m + 1}^2
     . \]
  In this case, we write $O_{2 m + 1} = O (V, Q_{2 m + 1})$.
\end{enumeratenumeric}

\subsection{Finite Euclidean Graphs}

Let $Q$ be a non-degenerate quadratic form on $V$. Then the finite Euclidean
graph $E_q (n, Q, a)$ is defined as the graph with vertex set $V$ and the edge
set

\begin{equation}
  E =\{(x, y) \in V \times V\, |\, x \neq y,\, Q (x - y) = a\}.
\end{equation}

In \cite{bannai 1}, Bannai, Shimabukuro and Tanaka showed that the finite Euclidean graphs $E_q(n,Q,a)$ are not always Ramanujan. Fortunately, they are always asymptotically Ramanujan. The following theorem summaries (in a rough form) the results from Sections~2-6 in \cite{bannai 1} and Section 3 in Kwok \cite{kwok}.

\begin{theorem} \label{euclidean graphs}
  Let $\rho$ be a primitive element of $\mathbbm{F}_q$.
  \begin{enumeratealpha}
    \item The graphs $E_q (2 m, Q^{\pm}_{2 m}, \rho^i)$ are regular of valency
    $k = q^{2 m - 1} \pm q^{m - 1}$ for $1 \leqslant i \leqslant q - 1$. Let
    $\lambda$ be any eigenvalue of the graph $E_q (2 m, Q^{\pm}_{2 m},
    \rho^i)$ with $\lambda \neq$ valency of the graph then
    \[ | \lambda | \leqslant 2 q^{(2 m - 1) / 2} . \]
    \item The graphs $E_q (2 m + 1, Q_{2 m + 1}, \rho^i)$ are regular of
    valency $k = q^{2 m} \pm q^m$ for $1 \leqslant i \leqslant q - 1$. Let
    $\lambda$ be any eigenvalue of the graph $E_q (2 m + 1, Q_{2 m + 1},
    \rho^i)$ with $\lambda \neq$ valency of the graph then
    \[ | \lambda | \leqslant 2 q^m . \]
  \end{enumeratealpha}
\end{theorem}

\subsection{Finite non-Euclidean Graphs}

In order to keep this paper concise, we will restrict our discussion to the finite
non-Euclidean graphs obtained from the action of the simple orthogonal group on the set of non-isotropic points.
Similar results hold for graphs obtained from the action of various Lie groups
on the set of non-isotropic points, lines and hyperplanes. We will address these results in a subsequent paper.

\subsubsection{Graphs obtained from the action of simple orthogonal group
$O_{2 m + 1} (q)$ ($q$ odd) on the set of non-isotropic points}

Let $V = \mathbbm{F}_q^{2 m + 1}$ be the $(2 m + 1)$-dimensional vector space
over the finite field $\mathbbm{F}_q$ ($q$ is an odd prime power). For each
element $x$ of $V$, we denote the $1$-dimensional subspace containing $x$ by
$[x]$. Let $\Theta, \Omega$ be the set of all square type and the set of all
non-square-type non-isotropic $1$-dimensional subspaces of $V$ with respect to
the quadratic form $Q_{2 m + 1}$, respectively. Then we have $| \Theta | =
(q^{2 m} - q^m) / 2$ and $| \Omega | = (q^{2 m} + q^m) / 2$. The simple
orthogonal group $O_{2 m + 1} (q)$ acts transitively on $\Theta$ and $\Omega$.

We define the graphs $H_q (O_{2 m + 1}, \Theta, i)$ (for $1 \leqslant i
\leqslant (q + 1) / 2$) as follows (let $E_i$ be the edge set of $H_q (O_{2 m
+ 1}, \Theta, i)$):
\begin{eqnarray*}
  ([x], [y]) \in E_1  & \Leftrightarrow & 
\left( \begin{array}{c}
					x\\
					y 
	\end{array} \right).S.
	\left( \begin{array}{c}
					x\\
					y 
	\end{array} \right)^{t} =
  \left( \begin{array}{cc}
					\nu & 1\\
					1 & \nu^{-1} 
	\end{array} \right),\\
  ([x], [y]) \in E_i & \Leftrightarrow & 
  \left( \begin{array}{c}
					x\\
					y 
	\end{array} \right).S.
	\left( \begin{array}{c}
					x\\
					y 
	\end{array} \right)^{t} = 
	\left( \begin{array}{cc}
					\nu & 1\\
					1 & \nu^{2i-3} 
	\end{array} \right), (2 \leqslant i \leqslant (q - 1) /
  2)\\
  ([x], [y]) \in E_{(q + 1) / 2} & \Leftrightarrow & 
  \left( \begin{array}{c}
					x\\
					y 
	\end{array} \right).S.
	\left( \begin{array}{c}
					x\\
					y 
	\end{array} \right)^{t} =
  \left( \begin{array}{cc}
					\nu & 0\\
					0 & \nu 
	\end{array} \right),
\end{eqnarray*}
where $\nu \in \mathbbm{F}_q$ is a primitive element of $\mathbbm{F}_q$, $A^t$
denotes the transpose of $A$ and $S$ is the matrix of the associated bilinear
form of $Q_{2 m + 1}$. Note that for $m = 1$ then we have the finite analog $H_q$
of the upper half plane.

We define the graph $H_q (O_{2 m + 1}, \Omega, i)$ (for $1 \leqslant i
\leqslant (q + 1) / 2)$ as follows (let $E_i$ be the edge set of $H_q (O_{2 m
+ 1}, \Omega, i)$):
\[\begin{array}{lll}
  ([x], [y]) \in E_1  & \Leftrightarrow & Q_{2 m + 1} (x + y) = 0,\\
  ([x], [y]) \in E_i & \Leftrightarrow & Q_{2 m + 1} (x + y) = 2 + 2 \nu^{- (i
  - 1)}, (2 \leqslant i \leqslant (q - 1) / 2)\\
  ([x], [y]) \in E_{(q + 1) / 2} & \Leftrightarrow & Q_{2 m + 1} (x + y) = 2,
\end{array}\]
where we assume $Q_{2 m + 1} (x) = 1$ for all $[x] \in \Omega$.

As in finite Euclidean case, the graphs obtained in this section are always asymptotically Ramanujan. The following theorem summaries the results from Sections 1 and 2 in \cite{bannai 2} and from Section 7 in \cite{bannai-hao-song}.

\begin{theorem} \label{non-euclidean graph 1}
  a) The graphs $H_q (O_{2 m + 1}, \Theta, i)$ $(1 \leqslant i \leqslant (q -
  1) / 2)$ are regular of valency $q^{2 m - 1} \pm q^{m - 1}$. The graph
  $H_q (O_{2 m + 1}, \Theta, (q + 1) / 2)$ is regular of valency $(q^{2 m - 1}
  \pm q^{m - 1})/2$. Let $\lambda$ be any eigenvalue of the graph $H_q
  (O_{2 m + 1}, \Theta, i)$ with $\lambda \neq$ valency of the graph then
  \[ | \lambda | \leqslant 2 q^{(2 m - 1) / 2} . \]
  b) The graphs $H_q (O_{2 m + 1}, \Omega, i)$ ($1 \leqslant i \leqslant (q -
  1) / 2$) are regular of valency $q^{2 m - 1} \pm q^{m - 1})$. The graph
  $H_q (O_{2 m + 1}, \Omega, (q + 1) / 2)$ is regular of valency $(q^{2 m - 1}
  \pm q^{m - 1})/2$. Let $\lambda$ be any eigenvalue of the graph $H_q
  (O_{2 m + 1}, \Omega, i)$ with $\lambda \neq$ valency of the graph then
  \[ | \lambda | \leqslant 2 q^{(2 m - 1) / 2} . \]
\end{theorem}

\subsubsection{Graphs obtained from the action of simple orthogonal group
$O_{2 m}^{\pm} (q)$ ($q$ odd) on the set of non-isotropic points}

Let $V = \mathbbm{F}_q^{2 m}$ be the $2 m$-dimensional vector space over the
finite field $\mathbbm{F}_q$ ($q$ is an odd prime power). For each element $x$
of $V$, we denote the $1$-dimensional subspace containing $x$ by $[x]$. Let
$\Omega_1, \Omega_2$ be the set of all square type and the set of all
non-square-type non-isotropic $1$-dimensional subspaces of $V$ with respect to
the quadratic form $Q^+_{2 m}$, respectively. Then we have $| \Omega_1 | = |
\Omega_2 | = (q^{2 m - 1} - q^{m - 1}) / 2$. The orthogonal group $O^+_{2 m}
(q)$ with respect to the quadratic from $Q^+_{2 m}$ over $\mathbbm{F}_q$ acts
on both $\Omega_1$ and $\Omega_2$ transitively. We define the graph $H_q (O_{2
m}^{+}, \Omega_1, i)$ (for $1 \leqslant i \leqslant (q + 1) / 2)$ as follows
(let $E_i$ be the edge set of $H_q (O_{2 m}^{+}, \Omega_1, i)$):

\[\begin{array}{lll}
  ([x], [y]) \in E_i & \Leftrightarrow & \left\langle x, y \right\rangle_{Q_{2
  m}^+} = 2^{- 1} \nu^i, (1 \leqslant i \leqslant (q - 1) / 2)\\
  ([x], [y]) \in E_{(q + 1) / 2} & \Leftrightarrow & \left\langle x, y
  \right\rangle_{Q_{2 m}^+} = 0,
\end{array}\]
where we assume $Q_{2 m}^+ (x) = 1$ for all $[x] \in \Omega$.

Let $\Theta_1, \Theta_2$ be the set of all square type and the set of all
non-square-type non-isotropic $1$-dimensional subspaces of $V$ with respect to
the quadratic form $Q^-_{2 m}$, respectively. Then we have $| \Theta_1 | = |
\Theta_2 | = (q^{2 m - 1} + q^{m - 1}) / 2$. The orthogonal group $O^-_{2 m}
(q)$ with respect to the quadratic from $Q^-_{2 m}$ over $\mathbbm{F}_q$ acts
on both $\Theta_1$ and $\Theta_2$ transitively. We define the graph $H_q (O_{2
m}^{-}, \Theta_1, i)$ (for $1 \leqslant i \leqslant (q + 1) / 2)$ as follows (let
$E_i$ be the edge set of $H_q (O_{2 m}^{-}, \Omega_1, i)$):

\[\begin{array}{lll}
  ([x], [y]) \in E_i & \Leftrightarrow & \left\langle x, y \right\rangle_{Q_{2
  m}^-} = 2^{- 1} \nu^i, (1 \leqslant i \leqslant (q - 1) / 2)\\
  ([x], [y]) \in E_{(q + 1) / 2} & \Leftrightarrow & \left\langle x, y
  \right\rangle_{Q_{2 m}^-} = 0,
\end{array}\]
where we assume $Q_{2 m}^- (x) = 1$ for all $[x] \in \Omega$.

The graphs obtained in this section are always asymptotically Ramanujan. The following theorem summaries the results from Sections 4 and 5 in \cite{bannai 2} and from Section 4 in \cite{bannai-hao-song}.

\begin{theorem} \label{non-euclidean graph 2}
  a) The graphs $H_q (O_{2 m}, \Theta_1, i)$ $(1 \leqslant i \leqslant (q - 1)
  / 2)$ are regular of valency $q^{2 m - 2} \pm q^{m - 1}$. The graph $H_q
  (O_{2 m}, \Theta, (q + 1) / 2)$ is regular of valency $(q^{2 m - 2} \pm
  q^{m - 1})/2$. Let $\lambda$ be any eigenvalue of the graph $H_q (O_{2 m},
  \Theta, i)$ with $\lambda \neq$ valency of the graph then
  \[ | \lambda | \leqslant 2 q^{(2 m - 2) / 2} . \]
  b) The graphs $H_q (O_{2 m}, \Omega_1, i)$ ($1 \leqslant i \leqslant (q - 1)
  / 2$) are regular of valency $q^{2 m - 2} \pm q^{m - 1}$. The graph $H_q
  (O_{2 m + 1}, \Omega, (q + 1) / 2)$ is regular of valency $(q^{2 m - 2} \pm q^{m - 1})/2$. Let $\lambda$ be any eigenvalue of the graph $H_q (O_{2 m},
  \Omega_1, i)$ with $\lambda \neq$ valency of the graph then
  \[ | \lambda | \leqslant 2 q^{(2 m - 2) / 2} . \]
\end{theorem}

\section{Explicit Tough Ramsey Graphs}
We call a graph $G = (V, E)$ $(n, d,
\lambda)$-regular if $G$ is a $d$-regular graph on $n$ vertices with the absolute value of each of its eigenvalues but the
largest one is at most $\lambda$. It is well-known that if $\lambda \ll d$ then a $(n,d,\lambda)$-regular graph behaves similarly as a random graph $G_{n,d/n}$. Presicely, we have the following result (see Corollary 9.2.5 and Corollary 9.2.6 in \cite{alon-spencer}).

\begin{theorem} \label{expander} (\cite{alon-spencer})
  Let $G$ be a $(n, d, \lambda)$-regular graph. 
  
  a) For every set of vertices $B$ and $C$ of $G$, we have
  \begin{equation} \label{f1}
    |e(B,C) - \frac{d}{n}|B||C| | \leqslant \lambda \sqrt{|B| |C|},
  \end{equation}
  where $e(B,C)$ is the number of edges in the induced subgraph of $G$ on $B$ (i.e. the number of ordered pairs $(u,v)$ where $u \in B, v \in C$ and $uv$ is an edge of $G$). 
  
  b) For every set of vertices $B$ of $G$, we have
  \begin{equation}\label{f2}
    |e (B) - \frac{d}{2 n} |B|^2 | \leqslant \frac{1}{2} \lambda |B|,
  \end{equation}
  where $e (B)$ is number of edges in the induced subgraph of $G$ on $B$.
\end{theorem}

Let $B, C$ be one of the maximum independent pairs of $G$, i.e. the ``bipartite'' subgraph induced on $(B,C)$ are empty and $|B||C|$ is maximum. Let $\alpha_2(G)$ denote the size $|B||C|$ of this pair. Then from (\ref{f1}), we have
\begin{equation} \label{independent pair}
\alpha_2(G) \leqslant \frac{\lambda^2n^2}{d^2}.
\end{equation}

Let $B$ be one of the maximum independent sets of $G$. Then from (\ref{f2}), we have
\begin{equation} \label{independent number}
\alpha(G) = |B| \leqslant \frac{n\lambda}{d},
\end{equation}
and
\begin{equation} \label{chromatic}
\chi(G) \geqslant \frac{|V(G)|}{\alpha(G)} \geqslant \frac{d}{\lambda}.
\end{equation}

The \textit{toughness} $t(G)$ of a graph $G$ is the largest real $t$ so that for every positive integer $x \geq 2$ one should delete at least $tx$ vertices from $G$ in order to get an induced subgraph of it with at least $x$ connected components. $G$ is $t$-tough if $t(G) \geq t$. This parameter was introduced by Chvat\'al in \cite{chvatal}. Chvat\'al also conjectures the following: \textit{there exists an absolute constant $t_0$ such that every $t_0$-tough graph is pancyclic.} This conjecture was disproved by Bauer, van den Heuvel and Schmeichel \cite{bauer} who constructed, for every real $t_0$, a $t_0$-tough triangle-free graph. They define a sequence of triangle-free graphs $H_1,H_2,H_3,\ldots$ with $|V(H_j)| = 2^{2j-1} (j+1)!$ and $t(H_j) \geq \sqrt{2j+4}/2$. To bound the toughness of a $(n,d,\lambda)$-regular graph, we have the following result which is due to Alon in \cite{alon1}.

\begin{theorem} \cite{alon1} \label{tough theorem} Let $G = (V,E)$ be an $(n,d,\lambda)$-graph. Then the toughness $t=t(G)$ of $G$ satisfies
\begin{equation} \label{tough}
t > \frac{1}{3} \left( \frac{d^2}{\lambda d + \lambda^2} - 1 \right).
\end{equation}
\end{theorem}

Let $G$ be any graph of the form $E_q (2 m, Q^{\pm}_{2 m}, a)$, $E_q (2 m + 1, Q_{2 m
+ 1}, a)$, $H_q (2 m + 1, \Theta, i)$, $H_q (2 m + 1, \Omega, i)$, $H_q (2 m,
\Omega_1, i)$ and $H_q (2 m, \Theta_1, i)$ for $a \neq 0 \in \mathbbm{F}_q$
and $1 \leqslant i \leqslant (q + 1) / 2$. Then from Theorems \ref{euclidean graphs}, \ref{non-euclidean graph 1} and \ref{non-euclidean graph 2}, the graph $G$ is $(c_1 q^n + O(q^{n/2}), c_2q^{n-1}+O(q^{(n-1)/2}),2q^{(n-1)/2})$-regular for some $n \geqslant 2$ and $c_1,c_2 \in \{\frac{1}{2},1\}$. By (\ref{independent number}), (\ref{chromatic}) and (\ref{tough}), we can show  that the finite Euclidean and non-Euclidean graphs have high chromatic number, small independent number and high tough number.

\begin{theorem} \label{bound}
Let $G$ be any graph of the form $E_q (2 m, Q^{\pm}_{2 m}, a)$, $E_q (2 m + 1, Q_{2 m
+ 1}, a)$, $H_q (2 m + 1, \Theta, i)$, $H_q (2 m + 1, \Omega, i)$, $H_q (2 m,
\Omega_1, i)$ and $H_q (2 m, \Theta_1, i)$ for $a \neq 0 \in \mathbbm{F}_q$
and $1 \leqslant i \leqslant (q + 1) / 2$. Suppose that $|V(G)| = cq^n + O(q^{(n-1)/2})$.
\begin{enumerate}
\item The independent number of $G$ is small: $\alpha(G) \leqslant (4+o(1)) |V(G)|^{(n+1)/2n}$.
\item The chromatic number of $G$ is high: $\chi(G) \geqslant |V(G)|^{(n-1)/2n}/(4+o(1))$.
\item The toughness of $G$ is at least $|V(G)|^{(n-1)/2n}/(12+o(1))$.
\end{enumerate}
\end{theorem}

In \cite{vinh}, the authors derived the following theorem using only elementary algebra. This theorem can also be derived from character tables of the association schemes of affine type (\cite{kwok}) and of finite orthogonal groups acting on the nonisotropic points (\cite{bannai-hao-song}). 

\begin{theorem} \label{triangle-free}
Among all finite Euclidean and non-Euclidean graphs, the only triangle-free graphs are
\begin{enumerate}
\item $E_q(2,Q^{-},a)$ where $3$ is square in $\mathbbm{F}_q$.
\item $E_q(2,Q^{+},a)$ where $3$ is nonsquare in $\mathbbm{F}_q$.
\item $H_q(3,Q,a)$ for at least one element $a \in \mathbbm{F}^{*}_q$.
\end{enumerate}
\end{theorem}

Theorems \ref{bound} and \ref{triangle-free} shows that the finite Euclidean $E_q(2,Q^{+},a)$, where $q$ is a prime of form $q=12k\pm5$ and $a \neq 0 \in \mathbbm{F}_q$, is an explicit triangle-free graph on $n_q=q^2$ vertices whose chromatic number exceeds $0.5n_q^{1/4}$. Therefore, this disproves the conjecture of Chavat\'al.  In addition, this graph is an explicit construction showing that $R(3,k) \geq \Omega(k^{4/3})$. 

The bounds obtained from Theorems \ref{bound} and \ref{triangle-free} match with the bounds obtained by code graphs in Theorem 3.1 in \cite{alon1}. These graphs are Caley graphs and their construction is based on some of the properties of certain Dual BCH error-correcting codes. For a positive integer $k$, let $F_k = GF(2^k)$ denote the finite field with $2^k$ elements. The elements of $F_k$ are represented by binary vectors of length $k$. If $a$ and $b$ are two such vectors, let $(a,b)$ denote their concatenation. Let $G_k$ be the graph whose vertices are all $n = 2^{2k}$ binary vectors of length $2k$, where two vectors $u$ and $v$ are adjacent if and only if there exists a non-zero $z \in F_k$ such that $u+v=(z,z^3)$ mod $2$ where $z^3$ is computed in the field $F_k$. Then $G_k$ is a $d_k=2^k-1$-regular graph on $n_k=2^{2k}$. Moreover, $G_k$ is triangle-free with independence number at most $2n^{3/4}$. Noga Alon gives a better bound $R(m,3) \geq \Omega(m^{3/2})$ in \cite{alon} by considering a graph with vertex set of all $n=2^{3k}$ binary vectors of length $3k$ (instead of all binary vectors of length $2k$). Suppose that $k$ is not divisible by $3$. Let $W_0$ be the set of all nonzero elements $\alpha \in F_k$ such that the leftmost bit in the binary representation of $\alpha^7$ is $0$, and let $W_1$ be the set of all nonzero elements $\alpha \in F_k$ for which the leftmost bit of $\alpha^7$ is $1$. Then $|W_0|=2^{k-1}-1$ and $|W_1|=2^{k-1}$. Let $G_n$ be the graph whose vertices are all $n=2^{3k}$ binary vectors of length $3k$, where two vectors $u$ and $v$ are adjacent if and only if there exist $w_0 \in W_0$ and $w_1 \in W_1$ such that $u+v=(w_0,w_0^3,w_0^5)+(w_1,w_1^3,w_1^5)$ where the powers are computed in the field $F_k$ and the addition is addition module $2$. Then $G_n$ is a $d_n=2^{k-1}(2^{k-1}-1)$-regular graph on $n=2^{3k}$ vertices. Moreover, $G_n$ is a triangle-free graph with independence number at most $(36+o(1))n^{2/3}$. The problem of finding better bounds for the chromatic number of finite Euclidean and non-Euclidean graphs on the plane and the upper half plane, respectively touches on an important question in graph theory: what is the greatest possible chromatic number for a triangle-free regular graph of order n? A possible approach is to consider the existence of sum-free varieties in high dimensional vector spaces over finite fields. We see that the varieties of degree two only give us triangle-free graphs on vector spaces of dimension two. We hope to address this problem for varieties of higher dimension in a subsequent paper. 

\section{Erd\"os distance problem}

\subsection{Proof of Theorem \ref{v1}}

Let $Q$ be any non-degenerate quadratic of $\mathbbm{F}_q^n$. Recall that the Euclidean graph $E_q (d,Q,a)$ was defined as the graph with vertex set
  $V$ and edge set
  \[ E =\{(x, y) \in V \times V | x \neq y, Q(x - y) = a\}. \]

\begin{lemma} \label{l1}
  Let $E \subset \mathbbm{F}_q^d$ such that $|E| \geqslant 3 q^{\frac{d +
  1}{2}} .$ Then $\Delta_Q (E) = \mathbbm{F}_q$.
\end{lemma}

\begin{proof}
  By Theorem \ref{euclidean graphs}, each graph $E_q (d,Q, a)$ is a $(q^d, q^{d - 1} \pm q^{\lfloor(d - 1) / 2\rfloor}, 2 q^{(d - 1) / 2})$-regular graph. By (\ref{independent number})  , for any
  $a \neq 0 \in \mathbbm{F}_q$, we have
  \begin{equation}
    \alpha (E_q (d, Q, a)) \leqslant \frac{2 q^{(3 d - 1) / 2}}{q^{d - 1} - q^{(d
    - 1) / 2}} \leqslant 3 q^{(d + 1) / 2} .
  \end{equation}
  Thus, if $|E| \geqslant 3 q^{\frac{d + 1}{2}}$ then $E$ is not an independent
  set of $E_q (d, Q, a)$, or equivalently there exist $x, y \in E$ such that $Q(x - y) = a$
  for any $a \in \mathbbm{F}_q$. This concludes the proof of the lemma.
\end{proof}

\begin{lemma} \label{l2}
  For any $0 < \varepsilon < 1 / 2$. Let $E \subset \mathbbm{F}_q$ such that
  $|E| \geqslant 3 q^{\frac{d}{2} + \varepsilon}$. Then
  \begin{equation}
    | \Delta_Q (E) | \geqslant q^{\frac{1}{2} + \varepsilon},
  \end{equation}
  for any $q \geqslant 6^{1 / (\varepsilon - 1 / 2)}$.
\end{lemma}

\begin{proof}
  By Theorem \ref{euclidean graphs}, each graph $E_q (d, Q, a)$ is a $(q^d, q^{d - 1} \pm 
  q^{\lfloor(d - 1) / 2\rfloor}, 2 q^{(d - 1 /) 2})$-regular graph. By (\ref{independent number}), the
  number of edges of $E_q (d, Q, a)$ in the induced subgraph on $E$ is at
  most
  \begin{equation} \label{15}
    e_{E_q (d, Q, a)} (E) \leqslant \frac{q^{d - 1} + q^{(d - 1) / 2}}{2 q^d}
    |E|^2 + q^{(d - 1) / 2} |E|.
  \end{equation}
  Suppose that $\# \Delta_Q (E) < q^{1 / 2 + \varepsilon}$. From (\ref{15}), we have
  \begin{eqnarray*}
    \binom{|E|}{2} & = & \sum_{a\in \Delta_Q(E)} e_{E_q(d,Q,a)}(E)\\
    & < & q^{1 / 2 + \varepsilon} \left\{ \frac{q^{d - 1} +
    q^{(d - 1) / 2}}{2 q^d} |E|^2 + q^{(d - 1) / 2} |E| \right\}\\
    & < & |E|q^{\varepsilon - \frac{1}{2}} \left\{ \left( \frac{1}{2} +
    \frac{1}{2} q^{- (d - 1) / 2} \right) |E| + q^{(d + 1) / 2} \right\},
  \end{eqnarray*}
  which implies that
  \begin{eqnarray*}
    q^{\frac{1}{2} - \varepsilon} (|E| - 1) & < & \left( 1 + q^{- (d - 1) / 2}
    \right) |E| + 2 q^{(d + 1) / 2}\\
    & \leqslant & (1 + q^{- 1 / 2} + \frac{2}{3} q^{\frac{1}{2} -
    \varepsilon}) |E|.
  \end{eqnarray*}
  Therefore, we have
  \begin{eqnarray*}
    q^{\frac{1}{2} - \varepsilon} & > & \left( \frac{1}{3} q^{\frac{1}{2} -
    \varepsilon} - 1 - q^{\frac{1}{2}} \right) |E|\\
    & \geqslant & \left( q^{\frac{1}{2} - \varepsilon} - 3 - 3q^{-
    \frac{1}{2}} \right) q^{\frac{d}{2} + \varepsilon},
  \end{eqnarray*}
  which is a contradiction if $q > 6^{1 / (1 / 2 - \varepsilon)}$. The lemma
  follows.
\end{proof}

Theorem \ref{v1} follows immediately from Lemma \ref{l1} and Lemma \ref{l2}.

\subsection{Proof of Theorem \ref{v2}}
For a fixed $a \in \mathbbm{F}_q$, the \textit{finite non-Euclidean graph} $V_q(\sigma,a)$ has vertices as the points in $H_q$ and edges between vertices $z,w$ if and only if $d(z,w) = a$. Except when $a=0$ or $a=4\sigma$, $V_q(\sigma,a)$ is a connected $(q+1)$-regular graph. When $a=0, 4\sigma$ then $V_q(\sigma,a)$ is disconnected, with one or two nodes, respectively, per connected component. As $a$ varies, we have $q-2$ $(q+1)$-regular graphs $V_q(\sigma,a)$. The question of whether these graphs are always nonisomorphic or not is still open.  See \cite{survey} for a survey of spectra of Laplacians of this graph. 

\begin{lemma} \label{l3}
  Let $E \subset H_q$ such that $|E| \geqslant 2 q^{3/2} .$ Then $| \Delta_H (E)| \geqslant q-1$.
\end{lemma}

\begin{proof}
  Each graph $V_q (\sigma, a)$ (with $a \neq 0, 4\sigma \in \mathbbm{F}_q$) is a $(q^2-q, q + 1, 2 q^{1/ 2})$-regular graph. By (\ref{independent number}), for any 
  $a \neq 0, 4\sigma \in \mathbbm{F}_q$, we have
  \begin{equation}
    \alpha (V_q (\sigma, a)) \leqslant \frac{2 (q^2-q)q^{1/2}}{q+1} \leqslant 2 q^{3 / 2} .
  \end{equation}
  Thus, $\#E \geqslant 2 q^{3/2}$ then $E$ is not an independent
  set of $V_q (\sigma, a)$ or equivalently, there exist $x, y \in E$ such that $d(x - y) = a$
  for any $a \in \mathbbm{F}_q - \{0, 4a\}$. This concludes the proof of the lemma.
\end{proof}

Note that $V_q(\sigma,4\sigma)$ is just a disjoint union of $(q^2-q)/2$ edges. So we can have a set $E \in H_q$ with $|E| =(q^2-q)/2$
and $\Delta_H(E) = \mathbbm{F}_q - \{4\sigma\}$.

\begin{lemma} \label{l4}
  For any $0 < \varepsilon < 1 / 2$. Let $E \subset \mathbbm{F}_q$ such that
  $|E| \geqslant 3 q^{\frac{d}{2} + \varepsilon}$. Then
  \begin{equation}
    | \Delta_H (E)| \geqslant q^{\frac{1}{2} + \varepsilon},
  \end{equation}
  for any $q \geqslant 9^{1 / (\varepsilon - 1 / 2)}$.
\end{lemma}

\begin{proof}
  For any $a \neq 0, 4\sigma \in \mathbbm{F}_q$, each graph $V_q (\sigma, a)$ is a $(q^2-q, q +1
  , 2 q^{1/2})$-regular graph. From Theorem \ref{expander}, the
  number of edges of $V_q (\sigma, a)$ in the induced subgraph on $E$ is at
  most
  \begin{equation} \label{18}
    e_{V_q (\sigma, a)} (E) \leqslant \frac{q+1}{2 (q^2-q)}
    |E|^2 + q^{1 / 2} |E|.
  \end{equation}
  Suppose that $|\Delta_H (E)| < q^{1 / 2 + \varepsilon}$. From (\ref{18}), we have
  \begin{eqnarray*}
    \binom{|E|}{2} &=& \sum_{a\in \Delta_H(E)} e_{V_q(\sigma,a)}(E)\\
    & < & q^{1 / 2 + \varepsilon} \left\{ \frac{q +
    1}{2 (q^2-q)} |E|^2 + q^{1 / 2} |E| \right\}\\
    & < & |E|q^{\varepsilon - \frac{1}{2}} \left\{ \left( \frac{1}{2} +
    \frac{1}{q-2} \right) |E| + q^{3 / 2} \right\},
  \end{eqnarray*}
  which implies that
  \begin{eqnarray*}
    q^{\frac{1}{2} - \varepsilon} (|E| - 1) & < & \left( 1 + \frac{2}{q-2}
    \right) |E| + 2 q^{3 / 2}\\
    & \leqslant & (1 + \frac{2}{q-2} + \frac{2}{3} q^{\frac{1}{2} -
    \varepsilon}) |E|.
  \end{eqnarray*}
  Therefore, we have
  \begin{eqnarray*}
    q^{\frac{1}{2} - \varepsilon} & > & \left( \frac{1}{3} q^{\frac{1}{2} -
    \varepsilon} - 1 - \frac{2}{q-2} \right) |E|\\
    & \geqslant & \left( q^{\frac{1}{2} - \varepsilon} - 3 - \frac{6}{q-2} \right) q^{1 + \varepsilon},
  \end{eqnarray*}
  which is a contradiction when $q > 9^{1 / (1 / 2 - \varepsilon)}$. The lemma
  follows.
\end{proof}

Theorem \ref{v2} follows immediately from Lemma \ref{l3} and Lemma \ref{l4}. Similar results hold for others non-Euclidean spaces defined in Section 2. We will discuss these results in a subsequent paper.

\subsection{Set of distances between two sets}

Now we will prove Theorem \ref{v3} and Theorem \ref{v4}. For any $a \neq 0 \in \mathbbm{F}_q$, by Theorem \ref{expander}, the number of edges of the graph $E_q(d,Q,a)$ in the induced ``bipartite'' subgraph on $(E,F)$ (two vertex parts are not necessary disjoint) is at most:
\begin{equation}\label{old} e_{E_q(d,Q,a)} \leqslant \frac{q^{d-1}+q^{(d-1)/2}}{q^d} |E| |F| + 2q^{(d-1)/2}\sqrt{|E||F|}. \end{equation}
Thus, we have 
\begin{eqnarray*}
 |E||F| & = & \sum_{a \in \Delta_Q(E,F)} {e_{E_q(d,Q,a)}}\\
        &\leqslant & \Delta_Q(E,F) \left( \frac{q^{d-1}+q^{(d-1)/2}}{q^d} |E| |F| + 2q^{(d-1)/2}\sqrt{|E||F|} \right),
\end{eqnarray*}
which implies that
\begin{equation}
\Delta_Q(E,F) \geqslant \frac{1}{\frac{1}{q} + \frac{1}{q^{(d+1)/2}} + \frac{2q^{(d-1)/2}}{\sqrt{|E||F|}}}.
\end{equation}

From the above inequality, we can easily derive the following analog of Lemma \ref{l2} for the distance set $\Delta_Q(E,F)$. 

\begin{lemma} \label{l5}
For any $0 < \epsilon < 1$. If $|E||F| \geqslant 9q^{(d-1)+\epsilon}$ then
\[ \Delta_Q(E,F) \geqslant \frac{\sqrt{|E||F|}}{3q^{(d-1)/2}} \geqslant q^{\epsilon/2}\]
for any $q \gg 1$.
\end{lemma}

By Theorem \ref{euclidean graphs}, each graph $E_q (d,Q, a)$ is a $(q^d, q^{d - 1} \pm q^{\lfloor(d - 1) / 2\rfloor}, 2 q^{(d - 1) / 2})$-regular graph. By (\ref{independent pair})  , for any
  $a \neq 0 \in \mathbbm{F}_q$, we have
  \begin{equation}
    \alpha_2 (E_q (d, Q, a)) \leqslant \left( \frac{2 q^{(3 d - 1) / 2}}{q^{d - 1} - q^{(d
    - 1) / 2}} \right)^2 \leqslant 9 q^{d + 1} .
  \end{equation}
  Thus, if $|E||F| \geqslant 9 q^{d + 1}$ then $E, F$ is not an independent
  pair of $E_q (d, Q, a)$ for any nonzero $a$. This implies that there exist $x \in E$ and $y \in F$ such that $Q(x,y) = a$
  for any $a \in \mathbbm{F}_q$. We have the following analog of Lemma \ref{l1}.
  
\begin{lemma} \label{l6}
  Let $E, F \subset \mathbbm{F}_q^d$ such that $|E||F| \geqslant 9 q^{d+1} .$ Then $\Delta_Q (E,F) = \mathbbm{F}_q$.
\end{lemma}

Theorem \ref{v3} is immediate from Lemma \ref{l5} and Lemma \ref{l6}. The proof of Theorem \ref{v4} is similar and is left for the readers. 
Note that the analog of Lemma \ref{l3} for the distance set $\Delta_H(E,F)$ is interesting in its own right.

\begin{lemma} \label{l7}
  Let $E, F \subset H_q$ such that $|E||F| \geqslant 9 q^{3} .$ Then $|\Delta_H (E,F)| \geq q-1$.
\end{lemma}

\textbf{Acknowledgments}\\
The author is very grateful to Dang Phuong Dung and Si Li for many useful discussions, helpful comments and endless encouragement.


\nocite{*}
\begin{thebibliography}{1}
\providecommand{\natexlab}[1]{#1}
\providecommand{\url}[1]{\texttt{#1}}
\expandafter\ifx\csname urlstyle\endcsname\relax
  \providecommand{\doi}[1]{doi: #1}\else
  \providecommand{\doi}{doi: \begingroup \urlstyle{rm}\Url}\fi


\bibitem{alon}
N. Alon, Explicit Ramsey graphs and orthonormal labellings, \textit{The Electronic Journal of Combinatorics} \textbf{1} (1994), R12, 8pp.

\bibitem{alon1}
N. Alon, Tough Ramsey graphs without short cycles, \textit{Journal of Algebraic Combinatorics}

\bibitem{alon-spencer}
N. Alon, J.H. Spencer, \textit{The Probabilisitc Method}, 2nd edition, Wiley-Interscience, 2000.

\bibitem{bannai-kwok-song}
E. Bannai, W.M. Kwok, S.-Y. Song, Ennola type dualities in the character tables of some association schemes, \textit{Mem. Fac. Sci. Kyushi Univ. Ser. A.} \textbf{44} (1990), 129-143.

\bibitem{bannai-hao-song}
E. Bannai, S. Hao, S.-Y. Song, Character tables of the association schemes of finite orthogonal groups acting on the nonisotropic points, \textit{Journal of Combinatorial Theory Series A} \textbf{54} (1990), 164-170.

\bibitem{bannai-hao-song-wei}
E. Bannai, S. Hao, S.-Y. Song, H. Wei, Character tables of certain association schemes coming from finite unitary and sympletic groups, \textit{Journal of Algebra} \textbf{144} (1991), 189-200.

\bibitem{bannai 2}
E. Bannai, O. Shimabukuro, H. Tanaka, Finite analogues of non-Euclidean spaces and Ramanujan graphs, \textit{European Journal of Combinatorics} \textbf{25} (2004), 243-259.

\bibitem{bannai 1}
E. Bannai, O. Shimabukuro, H. Tanaka, Finite Euclidean graphs and Ramanujan graphs, \textit{Discrete Mathematics} (to appear).

\bibitem{bauer}
D. Bauer, J. Vandenheuvel and E. Schmeichel, Toughness and Triangle-Free Graphs, \textit{Journal of Combinatorial Theory, Series B} \textbf{65} (2) (1995), 208-221.

\bibitem{bourgain}
J. Bourgain, Hausdorff dimension and distance sets, \textit{Israel Journal of Mathematics} \textbf{87} (1994), 1993-201.

\bibitem{bourgain-katz-tao}
J. Bourgain, N. Katz, T. Tao, A sum-product estimate in finite fields, and applications, \textit{Geom. Funct. Anal.} \textbf{14} (2004), 27-57.

\bibitem{chvatal}
V. Chv\'atal, Tough graphs and hamiltonian circuits, \textit{Discrete Mathematics} \textbf{5} (1973), 215-218.

\bibitem{hart-iosevich-solymosi}
D. Hart, A. Iosevich, J. Solymosi, Sum-product estimates in finite fields via Kloosterman sums, \textit{International Mathematics Research Notices} (to appear).

\bibitem{hart-iosevich}
D. Hart, A. Iosevich, Sums and products in finite fields: an integral geometric viewpoint, preprint.

\bibitem{iosevich-koh}
A. Iosevich, D. Koh, Erd\"os-Falconer distance problem, exponential sums, and Fourier analytic approach to incidence theorems in vector spaces over finite fields, preprint.

\bibitem{iosevich-rudnev}
A. Iosevich, M. Rudnev, Erd\"os distance problem in vector spaces over finite fields, \textit{Transactions of the American Mathematical Society} \textbf{359} (12) (2007), 6127-6142.

\bibitem{katz1}
N.M. Katz, Estimates for Soto-Andrade sums, \textit{J. Reine Angew. Math.} \textbf{438} (1993), 143-161.

\bibitem{katz2}
N.M. Katz, A note on exponential sums, \textit{Finite Fields Appl.} \textbf{1} (1995), 395-398.

\bibitem{kwok}
W.M. Kwok, Chracter table of association schemes of affine type, \textit{European Journal of Combinatorics}, \textbf{13} (1992), 167-185.

\bibitem{li}
W.-C.W. Li, \textit{Number Theory with Applications}, World Scientific, River Edge, NJ, 1996.

\bibitem
{before} A. Medrano, P. Myers, H.M. Stark and A. Terras, Finite analogues of Euclidean space,
\textit{Journal of Computational and Applied Mathematics}, \textbf{68} (1996), 221-238.

\bibitem{shparlinski}
I. E. Shparlinski, On the set of distances between two sets over finite fields, \textit{International Journal of Mathematics and Mathematical Sciences} \textbf{2006}, Article ID 59482, 1-5.

\bibitem{tao1}
T. Tao, Finite field analogues of Erd\"os, Falconer, and Furstenberg problems, preprint.

\bibitem{survey}
A. Terras, Survey of Spectra of Laplacians on Finite Symmetric Spaces, \textit{Experimental Mathematics}, (1996).

\bibitem{vinh1}
L.A. Vinh, Some coloring problems for unit-quadrance graphs,
\textit{The proceedings of Australian Workshop on Combinatorial Algorithms}, 2006, 361-367.

\bibitem{vinh}
L.A. Vinh and D.P.Dung, Explicit tough Ramsey graphs, preprint.

\bibitem{van}
V. Vu, Sum-Product estimates via directed expanders, (preprint), (2007).

\end{thebibliography}
\end{document}